\documentclass{article}
\usepackage{amsmath,amsfonts,amssymb,bm,hyperref,amsthm}
\numberwithin{equation}{section}
\newtheorem{theorem}{Theorem}[section]
\newtheorem{lemma}[theorem]{Lemma}

\theoremstyle{definition}

\newtheorem{remark}[theorem]{Remark}
\newtheorem{definition}[theorem]{Definition}

\usepackage{authblk}

\begin{document}

\title{Analysis as a source of geometry: a non-geometric representation of the Dirac equation}

\author{Yan-Long Fang\thanks{yan.fang.12@ucl.ac.uk}\ }
\author{Dmitri Vassiliev\thanks{D.Vassiliev@ucl.ac.uk, \url{http://www.homepages.ucl.ac.uk/\~ucahdva/} \ \ \ Supported by \\ EPSRC grant EP/M000079/1}}
\affil{Department of Mathematics,
University College London,\\ Gower Street, London WC1E 6BT, UK}

\maketitle

\begin{abstract}
Consider a formally self-adjoint first order linear differential
operator acting on pairs (2-columns) of complex-valued scalar fields
over a 4-manifold without boundary. We examine the geometric content of such an
operator and show that it implicitly contains a Lorentzian
metric, Pauli matrices, connection coefficients for spinor
fields and an electromagnetic covector potential. This
observation allows us to give a simple representation of the
massive Dirac equation as a system of four scalar equations
involving an arbitrary two-by-two matrix operator as above and its
adjugate. The point of the paper is that in order to write
down the Dirac equation in the physically meaningful
4-dimensional hyperbolic setting one does not need any geometric
constructs. All the geometry required is contained in a single
analytic object --- an abstract formally self-adjoint first
order linear differential operator acting on pairs of
complex-valued scalar fields.
\end{abstract}


\section{Introduction}
\label{Introduction}

The paper is an attempt at developing a relativistic field theory based on
the concepts from the analysis of partial differential equations as opposed
to geometric concepts. The long-term goal is to recast quantum electrodynamics
in curved spacetime
in such ``non-geometric'' terms. The potential advantage of formulating a field theory
in ``analytic'' terms is that there might be a chance of describing
the interaction of different physical fields in a more consistent, and,
hopefully, non-perturbative manner.

The current paper deals with the Dirac equation in curved spacetime,
with the electromagnetic field appearing as a prescribed external covector potential.
We expect to treat the Maxwell system in a separate paper.

Let $M$ be a 4-manifold without boundary and let $m$ be the electron mass.

The traditional way of writing the massive Dirac equation
is as follows. We equip our manifold $M$ with a prescribed Lorentzian metric
and a prescribed electromagnetic covector potential, and write the Dirac
equation using the rules of spinor calculus, see
Appendix~\ref{Dirac equation in its traditional form}.
In the
process of doing this one may encounter topological
obstructions: not every 4-manifold admits a Lorentzian metric
and, even if it admits one, it may still not admit a spin
structure.

We give now an analytic representation of the massive
Dirac equation which, for parallelizable manifolds,
turns out to be equivalent to the traditional geometric representation.

For the sake of clarity, prior to describing our analytic construction
let us explain why we will not encounter topological obstructions
related to the second Stiefel--Whitney class. We will work with
operators satisfying the non-degeneracy condition (\ref{definition of non-degeneracy})
which is very natural from the analytic point of view as it is a generalisation
(weaker version) of the standard ellipticity condition (\ref{definition of ellipticity}).
It turns out that the imposition of the non-degeneracy condition
(\ref{definition of non-degeneracy}) has far reaching geometric consequences:
it implies that our manifold $M$ is parallelizable. Thus, in our
construction we deal only with parallelizable manifolds, but we do not
state the parallelizability condition explicitly because it is automatically
encoded in the analytic non-degeneracy condition
(\ref{definition of non-degeneracy}).

We assume that our 4-manifold $M$ is equipped with a
prescribed positive density $\rho$ which allows us to define
an inner product on columns of complex-valued scalar fields,
see formula (\ref{definition of sesquilinear and bilinear forms}),
and, consequently, the concept of formal self-adjointness,
see formula (\ref{definition of adjoint and transpose}).

Let $L$ be a first order linear
differential operator acting on 2-columns of complex-valued
scalar fields over $M$.
The standard invariant analytic way of describing this operator
is by means of its principal symbol $L_\mathrm{prin}(x,p)$
and subprincipal symbol $L_\mathrm{sub}(x)$,
see Appendix~\ref{Basic notions from the analysis of PDEs} for details.
Here $x=(x^1,x^2,x^3,x^4)$ are local coordinates on $M$
and $p=(p_1,p_2,p_3,p_4)$ is the dual variable (momentum).
It is known that $L_\mathrm{prin}$ and $L_\mathrm{sub}$
are invariantly defined
$2\times2$
matrix-functions on $T^*M$ and $M$
respectively and that these matrix-functions completely determine the
first order differential operator $L$.

Further on we assume that our differential operator $L$
is formally self-adjoint  and satisfies the
non-degeneracy condition (\ref{definition of non-degeneracy}).

We now take an arbitrary matrix-function
\begin{equation}
\label{GL2C matrix-function Q}
Q:M\to\mathrm{GL}(2,\mathbb{C})
\end{equation}
and consider the transformation of our differential operator
\begin{equation}
\label{GL2C transformation of the operator}
L\mapsto Q^*LQ.
\end{equation}

The motivation for looking at such transformations
is as follows. Let us write down the action (variational functional)
associated with our operator,
$\int_M v^*(Lv)\,\rho\,dx\,$,
and let us perform an invertible linear transformation
\begin{equation*}
\label{motivation 1}
v\mapsto Qv
\end{equation*}
in the vector space
$V:=\{v:M\to\mathbb{C}^2\}$
of 2-columns of complex-valued scalar fields.
Then the action transforms as
\begin{equation*}
\label{motivation 2}
\int_M v^*(Lv)\,\rho\,dx
\,\mapsto\,
\int_M v^*(Q^*LQv)\,\rho\,dx\,.
\end{equation*}
We see that the transformation
(\ref{GL2C transformation of the operator})
of our differential operator describes the transformation of the
integrand in the formula for the action.
We choose to interpret
(\ref{GL2C transformation of the operator}) as a gauge transformation.

The transformation (\ref{GL2C transformation of the operator})
of the differential operator $L$ induces the following transformations
of its principal and subprincipal symbols:
\begin{equation}
\label{GL2C transformation of the principal symbol}
L_\mathrm{prin}\mapsto Q^*L_\mathrm{prin}Q,
\end{equation}
\begin{equation}
\label{GL2C transformation of the subprincipal symbol}
L_\mathrm{sub}\mapsto
Q^*L_\mathrm{sub}Q
+\frac i2
\left(
Q^*_{x^\alpha}(L_\mathrm{prin})_{p_\alpha}Q
-
Q^*(L_\mathrm{prin})_{p_\alpha}Q_{x^\alpha}
\right),
\end{equation}
where the subscripts indicate partial derivatives.
Here we made use of formula (9.3) from \cite{jst_part_a}.

Comparing formulae
(\ref{GL2C transformation of the principal symbol})
and
(\ref{GL2C transformation of the subprincipal symbol})
we see that, unlike the principal symbol, the subprincipal
symbol does not transform in a covariant fashion due to
the appearance of terms with the gradient of the
matrix-function $Q(x)$. In order to identify the sources of this
non-covariance we observe that
any matrix-function (\ref{GL2C matrix-function Q})
can be written as a product of three terms:
a complex matrix-function of determinant one,
a positive scalar function and
a complex scalar function of modulus one.
Hence, we examine the three gauge-theoretic actions separately.

Take an arbitrary scalar function
\begin{equation}
\label{function psi}
\psi:M\to\mathbb{R}
\end{equation}
and consider the transformation of our differential operator
\begin{equation}
\label{psi transformation of the operator}
L\mapsto e^\psi Le^\psi.
\end{equation}
The transformation
(\ref{psi transformation of the operator}) is a special case of the transformation
(\ref{GL2C transformation of the operator}) with $Q=e^\psi I$,
where $I$ is the $2\times2$ identity matrix.
Substituting this $Q$ into formula
(\ref{GL2C transformation of the subprincipal symbol}),
we get
\begin{equation}
\label{psi transformation of the subprincipal symbol}
L_\mathrm{sub}\mapsto
e^{2\psi}L_\mathrm{sub},
\end{equation}
so the subprincipal symbol transforms in a covariant fashion.

Now take an arbitrary scalar function
\begin{equation}
\label{function phi}
\phi:M\to\mathbb{R}
\end{equation}
and consider the transformation of our differential operator
\begin{equation}
\label{phi transformation of the operator}
L\mapsto e^{-i\phi}Le^{i\phi}.
\end{equation}
The transformation
(\ref{phi transformation of the operator}) is a special case of the transformation
(\ref{GL2C transformation of the operator}) with $Q=e^{i\phi}I$.
Substituting this $Q$ into formula
(\ref{GL2C transformation of the subprincipal symbol}),
we get
\begin{equation}
\label{phi transformation of the subprincipal symbol}
L_\mathrm{sub}(x)\mapsto
L_\mathrm{sub}(x)+L_\mathrm{prin}(x,(\operatorname{grad}\phi)(x)),
\end{equation}
so the subprincipal symbol does not transform in a covariant fashion.
We do not take any action with regards to the non-covariance
of (\ref{phi transformation of the subprincipal symbol}).

Finally, take an arbitrary matrix-function
\begin{equation}
\label{SL2C matrix-function R}
R:M\to\mathrm{SL}(2,\mathbb{C})
\end{equation}
and consider the transformation of our differential operator
\begin{equation}
\label{SL2C transformation of the operator}
L\mapsto R^*LR.
\end{equation}
Of course, the transformation
(\ref{SL2C transformation of the operator}) is a special case of the transformation
(\ref{GL2C transformation of the operator}): we are looking at the
case when $\det Q(x)=1$.
It turns out that it is possible to overcome the resulting non-covariance
in (\ref{GL2C transformation of the subprincipal symbol}) by introducing
the \emph{covariant subprincipal symbol} $\,L_\mathrm{csub}(x)\,$
in accordance with formula
\begin{equation}
\label{definition of covariant subprincipal symbol}
L_\mathrm{csub}:=
L_\mathrm{sub}-f(L_\mathrm{prin}),
\end{equation}
where $f$ is a function (more precisely, a nonlinear differential operator)
mapping a $2\times2$ non-degenerate Hermitian principal symbol
$L_\mathrm{prin}(x,p)$
to a $2\times2$ Hermitian matrix-function
$(f(L_\mathrm{prin}))(x)$.
The function $f$ is chosen from the condition that
the transformation (\ref{SL2C transformation of the operator})
of the differential operator induces the transformation
\begin{equation}
\label{SL2C transformation of the covariant subprincipal symbol}
L_\mathrm{csub}\mapsto
R^*L_\mathrm{csub}R
\end{equation}
of its covariant subprincipal symbol and the condition
\begin{equation}
\label{homogeneity condition for f}
f(e^{2\psi}L_\mathrm{prin})=e^{2\psi}f(L_\mathrm{prin}),
\end{equation}
where $\psi$ is an arbitrary scalar function (\ref{function psi}).

The existence of a function $f$ satisfying conditions
(\ref{SL2C transformation of the covariant subprincipal symbol})
and
(\ref{homogeneity condition for f}) is a nontrivial fact,
a feature specific to a system of two equations
in dimension four. The explicit formula for the function $f$ is formula
(\ref{formula for f}).

Let us summarise the results of our gauge-theoretic analysis.
\begin{itemize}
\item
Our first order differential operator $L$ is completely
determined by its principal symbol $L_\mathrm{prin}(x,p)$
and covariant subprincipal symbol $L_\mathrm{csub}(x)$.
\item
The transformation
(\ref{GL2C transformation of the operator})
of the differential operator induces the transformation
(\ref{GL2C transformation of the principal symbol})
of its principal symbol.
\item
Transformations
(\ref{psi transformation of the operator}),
(\ref{phi transformation of the operator})
and
(\ref{SL2C transformation of the operator})
of the differential operator induce transformations
\begin{equation}
\label{psi transformation of the covariant subprincipal symbol}
L_\mathrm{csub}\mapsto
e^{2\psi}L_\mathrm{csub}\,,
\end{equation}
\begin{equation}
\label{phi transformation of the covariant subprincipal symbol}
L_\mathrm{csub}(x)\mapsto
L_\mathrm{csub}(x)
+L_\mathrm{prin}(x,(\operatorname{grad}\phi)(x))
\end{equation}
and (\ref{SL2C transformation of the covariant subprincipal symbol})
of its covariant subprincipal symbol.
\end{itemize}

We use the notation
\begin{equation}
\label{symbol Op}
L=\operatorname{Op}(L_\mathrm{prin},L_\mathrm{csub})
\end{equation}
to express the fact that our operator is completely
determined by its principal symbol
and covariant subprincipal symbol.
The differential operator $L$ can be written down explicitly,
in local coordinates,
via the principal symbol $L_\mathrm{prin}$
and covariant subprincipal symbol $L_\mathrm{csub}$
in accordance with formula (\ref{symbol Op explicit}),
so formula (\ref{symbol Op}) is shorthand for (\ref{symbol Op explicit}).
We call (\ref{symbol Op}) the \emph{covariant representation}
of the differential operator $L$.

Recall now a definition from elementary linear algebra.
The \emph{adjugate} of a $2\times2$
matrix is defined as
\begin{equation}
\label{definition of adjugation}
P=\begin{pmatrix}a&b\\ c&d\end{pmatrix}
\mapsto
\begin{pmatrix}d&-b\\-c&a\end{pmatrix}
=:\operatorname{adj}P.
\end{equation}

Using the covariant representation (\ref{symbol Op})
and matrix adjugation (\ref{definition of adjugation})
we can define the
adjugate of the differential operator $L$ as
\begin{equation}
\label{definition of adjugate operator}
\operatorname{Adj}L:=\operatorname{Op}
(\operatorname{adj}L_\mathrm{prin},
\operatorname{adj}L_\mathrm{csub}).
\end{equation}

Note that in the case when the principal symbol does not depend
on the position variable $x$
(this corresponds to Minkowski spacetime,
which is the case most important for applications)
the definition of the adjugate differential operator
simplifies. In this case the subprincipal symbol coincides with the
covariant subprincipal symbol and
one can treat the differential operator $L$ as if
it were a matrix: formula (\ref{definition of adjugate operator}) becomes
\begin{equation}
\label{definition of adjugate operator simplified}
L=\begin{pmatrix}L_{11}&L_{12}\\ L_{21}&L_{22}\end{pmatrix}
\mapsto
\begin{pmatrix}L_{22}&-L_{12}\\-L_{21}&L_{11}\end{pmatrix}
=\operatorname{Adj}L.
\end{equation}

We define the Dirac operator as the differential operator
\begin{equation}
\label{analytic definition of the Dirac operator}
D:=
\begin{pmatrix}
L&mI\\
mI&\operatorname{Adj}L
\end{pmatrix}
\end{equation}
acting on 4-columns $v$ of complex-valued scalar fields.
Here $I$ is the $2\times2$ identity matrix.
We claim
that the system of four scalar equations
\begin{equation}
\label{analytic statement of the Dirac equation}
Dv=0
\end{equation}
is equivalent to the Dirac equation in its traditional geometric formulation.

Examination of formula
(\ref{analytic definition of the Dirac operator})
raises the following questions.
\begin{itemize}
\item
Where is the Lorentzian metric?
\item
Why don't we encounter topological obstructions?
\item
Where are the Pauli matrices?
\item
Where are the spinors?
\item
Where are the connection coefficients for spinor fields?
\item
Where is the electromagnetic covector potential?
\item
Where is Lorentz invariance?
\end{itemize}

These questions will be answered in
Sections~\ref{Lorentzian metric}--\ref{Lorentz invariance}.
In Section~\ref{Main result}
we will collect together all the
formulae from
Sections~\ref{Lorentzian metric}--\ref{Lorentz invariance}
and show, by direct substitution, that our
equation (\ref{analytic statement of the Dirac equation})
is indeed the Dirac equation (\ref{Dirac equation traditional}).
This fact will be presented in the form of Theorem~\ref{main theorem},
the main result of our paper.

\section{Lorentzian metric}
\label{Lorentzian metric}

Observe that the determinant of the principal symbol is a quadratic form
in the dual variable (momentum) $p$\,:
\begin{equation}
\label{definition of metric}
\det L_\mathrm{prin}(x,p)=-g^{\alpha\beta}(x)\,p_\alpha p_\beta\,.
\end{equation}
We interpret the real coefficients $g^{\alpha\beta}(x)=g^{\beta\alpha}(x)$,
$\alpha,\beta=1,2,3,4$, appearing in formula (\ref{definition of metric})
as components of a (contravariant) metric tensor.

\begin{lemma}
\label{Lemma about Lorentzian metric}
Our metric is Lorentzian, i.e.~it has
three positive eigenvalues
and one negative eigenvalue.
\end{lemma}

\emph{Proof\ }
Decomposing $L_\mathrm{prin}(x,p)$ with respect to the standard basis
\begin{equation}
\label{standard basis}
s^1=
\begin{pmatrix}
0&1\\
1&0
\end{pmatrix},
\quad
s^2=
\begin{pmatrix}
0&-i\\
i&0
\end{pmatrix},
\quad
s^3=
\begin{pmatrix}
1&0\\
0&-1
\end{pmatrix},
\quad
s^4=
\begin{pmatrix}
1&0\\
0&1
\end{pmatrix}
\end{equation}
in the real vector space of $2\times2$ Hermitian matrices, we get
\begin{equation}
\label{principal symbol via frame}
L_\mathrm{prin}(x,p)=s^j c_j(x,p),
\end{equation}
where the repeated index $j$ indicates summation over $j=1,2,3,4$
and the $c_j(x,p)$ are some real-valued functions on $T^*M$.
Each coefficient $c_j(x,p)$ is linear in $p$, so
\begin{equation}
\label{principal symbol via frame 1}
c_j(x,p)=e_j{}^\alpha(x)\,p_\alpha\,,
\end{equation}
where the repeated index $\alpha$ indicates summation over $\alpha=1,2,3,4$
and $e_j$ is some real-valued vector field with components $e_j{}^\alpha(x)$.
The quartet of real-valued vector fields $e_j$, $j=1,2,3,4$,
is called the \emph{frame}.
Note that the non-degeneracy condition
(\ref{definition of non-degeneracy})
ensures that the vector fields $e_j$ are linearly independent
at every point of our manifold $M$.

Substituting
(\ref{standard basis})
and
(\ref{principal symbol via frame 1})
into
(\ref{principal symbol via frame}), we get
\begin{equation}
\label{principal symbol via frame 2}
L_\mathrm{prin}(x,p)=s^j e_j{}^\alpha(x)\,p_\alpha=
\begin{pmatrix}
e_4{}^\alpha p_\alpha+e_3{}^\alpha p_\alpha&
e_1{}^\alpha p_\alpha-ie_2{}^\alpha p_\alpha\\
e_1{}^\alpha p_\alpha+ie_2{}^\alpha p_\alpha
&e_4{}^\alpha p_\alpha-e_3{}^\alpha p_\alpha
\end{pmatrix}.
\end{equation}
Calculating the determinant of  (\ref{principal symbol via frame 2})
and substituting the result into the LHS of (\ref{definition of metric}),
we get
$
g^{\alpha\beta}\,p_\alpha p_\beta
=
(e_1{}^\alpha p_\alpha)^2
+
(e_2{}^\alpha p_\alpha)^2
+
(e_3{}^\alpha p_\alpha)^2
-
(e_4{}^\alpha p_\alpha)^2
$.~$\square$

\

The proof of Lemma~\ref{Lemma about Lorentzian metric}
explains why we do not encounter topological obstructions:
condition (\ref{definition of non-degeneracy})
implies that our manifold is parallelizable.

It is also easy to see that our frame defined in accordance with
formula
(\ref{principal symbol via frame 2})
is orthonormal with respect to the metric
(\ref{definition of metric}):
\begin{equation}
\label{orthonormality of the frame}
g_{\alpha\beta}\,e_j{}^\alpha e_k{}^\beta=
\begin{cases}
0\quad\text{if}\quad j\ne k,
\\
1\quad\text{if}\quad j=k\ne4,
\\
-1\quad\text{if}\quad j=k=4.
\end{cases}
\end{equation}

\section{Geometric meaning of our transformations}

In Section~\ref{Introduction} we defined four
transformations of a formally self-adjoint $2\times2$
first order linear differential operator:
\begin{itemize}
\item
conjugation (\ref{psi transformation of the operator})
by a positive scalar function,
\item
conjugation (\ref{phi transformation of the operator})
by a complex scalar function of modulus one,
\item
conjugation (\ref{SL2C transformation of the operator})
by an $\mathrm{SL}(2,\mathbb{C})$-valued matrix-function
and
\item
adjugation (\ref{definition of adjugate operator}).
\end{itemize}
In this section we establish the geometric meaning of the transformations
(\ref{psi transformation of the operator}),
(\ref{SL2C transformation of the operator})
and
(\ref{definition of adjugate operator}).
We do this by looking at the resulting transformations of the principal symbol.

We choose to examine the three transformations listed above in reverse order:
first (\ref{definition of adjugate operator}),
then \ref{SL2C transformation of the operator})
and, finally, (\ref{psi transformation of the operator}).

We know that $L_\mathrm{prin}$
can be written in terms of the standard basis (\ref{standard basis})
and frame $e_j$ as  (\ref{principal symbol via frame 2}).
Similarly, $\operatorname{adj}L_\mathrm{prin}$
can be written as
\begin{equation}
\label{principal symbol via frame adjugate}
\operatorname{adj}L_\mathrm{prin}(x,p)
=s^j\tilde e_j{}^\alpha(x)\,p_\alpha=
\begin{pmatrix}
\tilde e_4{}^\alpha p_\alpha+\tilde e_3{}^\alpha p_\alpha&
\tilde e_1{}^\alpha p_\alpha-i\tilde e_2{}^\alpha p_\alpha\\
\tilde e_1{}^\alpha p_\alpha+i\tilde e_2{}^\alpha p_\alpha
&\tilde e_4{}^\alpha p_\alpha-\tilde e_3{}^\alpha p_\alpha
\end{pmatrix},
\end{equation}
where $\tilde e_j$ is another frame.
Examination of formulae
(\ref{definition of adjugation}),
(\ref{principal symbol via frame 2})
and
(\ref{principal symbol via frame adjugate})
shows that the two frames,
$e_j$ and $\tilde e_j$,
differ by spatial inversion:
\begin{equation}
\label{spatial inversion}
e_j\mapsto-e_j,\quad j=1,2,3,
\qquad
e_4\mapsto e_4.
\end{equation}

The transformation (\ref{SL2C transformation of the operator})
of the differential operator induces the following transformation
of its principal symbol:
\begin{equation}
\label{SL2C transformation of the principal symbol}
L_\mathrm{prin}\mapsto R^*L_\mathrm{prin}R.
\end{equation}
If we recast the transformation
(\ref{SL2C transformation of the principal symbol})
in terms of the frame $e_j$
(see formula (\ref{principal symbol via frame 2})),
we will see that we are looking at a linear transformation of the frame,
\begin{equation}
\label{Lorentz transformation of the frame}
e_j\mapsto\Lambda_j{}^ke_k\,,
\end{equation}
with some real-valued coefficients
$\Lambda_j{}^k(x)$. The transformation of the principal symbol
(\ref{SL2C transformation of the principal symbol})
preserves the Lorentzian metric (\ref{definition of metric}),
so the linear transformation of the frame (\ref{Lorentz transformation of the frame})
is a Lorentz transformation.

Of course, the transformation (\ref{spatial inversion})
is also a Lorentz transformation and it can be written in the
form (\ref{Lorentz transformation of the frame}) with
$\Lambda_j{}^k=\operatorname{diag}(-1,-1,-1,+1)$.
The difference between the two Lorentz transformations
is that in the case of
adjugation (\ref{definition of adjugate operator})
we get $\det \Lambda_j{}^k=-1$,
whereas in the case of
conjugation (\ref{SL2C transformation of the operator})
by an $\mathrm{SL}(2,\mathbb{C})$-valued matrix-function
we get $\det \Lambda_j{}^k=+1$.

Finally, let us establish the geometric meaning of
conjugation (\ref{psi transformation of the operator})
by a positive scalar function.
The transformation (\ref{psi transformation of the operator})
of the differential operator induces the following transformation
of its principal symbol:
\begin{equation}
\label{psi transformation of the principal symbol}
L_\mathrm{prin}\mapsto e^{2\psi}L_\mathrm{prin}.
\end{equation}
Comparing formulae
(\ref{definition of metric})
and
(\ref{psi transformation of the principal symbol})
we see that we are looking at a conformal scaling of the metric,
\begin{equation}
\label{conformal scaling of the metric}
g^{\alpha\beta}\mapsto e^{4\psi}g^{\alpha\beta}.
\end{equation}

\begin{remark}
We did not examine in this section the geometric meaning of the transformation
(\ref{phi transformation of the operator}).
We did not do it because this transformation
does not affect the principal symbol: one has to look at
the subprincipal symbol to understand the geometric meaning
of the transformation
(\ref{phi transformation of the operator}).
We will do this later, in
Section~\ref{Electromagnetic covector potential}:
see formula (\ref{phi transformation of A}).
\end{remark}

\section{Pauli matrices}
\label{Pauli matrices and spinors}

The principal symbol
$L_\mathrm{prin}(x,p)$ of our operator $L$
is linear in the dual variable $p$, so it can be written as
\begin{equation}
\label{definition of Pauli matrices}
L_\mathrm{prin}(x,p)=\sigma^\alpha(x)\,p_\alpha\,.
\end{equation}
The four matrix-functions $\sigma^\alpha(x)$, $\alpha=1,2,3,4$,
appearing in (\ref{definition of Pauli matrices}) are, by definition, our Pauli matrices.

The adjugate of the principal symbol can be written as
\begin{equation}
\label{definition of adjugate Pauli matrices}
\operatorname{adj}L_\mathrm{prin}(x,p)=\tilde\sigma^\alpha(x)\,p_\alpha\,.
\end{equation}
The matrices
$\tilde\sigma^\alpha(x)$, $\alpha=1,2,3,4$,
appearing in formula (\ref{definition of adjugate Pauli matrices})
are the adjugates of those from (\ref{definition of Pauli matrices})

We have
\begin{equation}
\label{Pauli extra formula}
[L_\mathrm{prin}(x,p)]
[\operatorname{adj}L_\mathrm{prin}(x,p)]
\!=\!
[\operatorname{adj}L_\mathrm{prin}(x,p)]
[L_\mathrm{prin}(x,p)]
\!=\!-Ig^{\alpha\beta}p_\alpha p_\beta,
\end{equation}
where $I$ is the $2\times2$ identity matrix
and $g^{\alpha\beta}$ is the metric from formula
(\ref{definition of metric}).
Formula (\ref{Pauli extra formula}) implies
\begin{equation*}
\label{defining identity for Pauli matrices 2}
[L_\mathrm{prin}(x,p)]
[\operatorname{adj}L_\mathrm{prin}(x,q)]
+
[L_\mathrm{prin}(x,q)]
[\operatorname{adj}L_\mathrm{prin}(x,p)]
=-2Ig^{\alpha\beta}p_\alpha q_\beta\,,
\end{equation*}
\begin{equation*}
\label{defining identity for Pauli matrices 3}
[\operatorname{adj}L_\mathrm{prin}(x,p)]
[L_\mathrm{prin}(x,q)]
+
[\operatorname{adj}L_\mathrm{prin}(x,q)]
[L_\mathrm{prin}(x,p)]
=-2Ig^{\alpha\beta}p_\alpha q_\beta\,.
\end{equation*}
Substituting
(\ref{definition of Pauli matrices})
and
(\ref{definition of adjugate Pauli matrices})
into the above formulae we arrive at
(\ref{defining identity for Pauli matrices})
and
(\ref{defining identity for Pauli matrices 1}).
This means that our matrices $\sigma^\alpha(x)$
defined in accordance with formula
(\ref{definition of Pauli matrices})
satisfy the abstract definition of Pauli matrices,
Definition~\ref{abstract definition of Pauli matrices}.

\section{Covariant subprincipal symbol}
\label{Covariant subprincipal symbol}

Recall that we defined the covariant subprincipal symbol
$L_\mathrm{csub}(x)$
in accordance with formula
(\ref{definition of covariant subprincipal symbol}).
We need now to determine the function $f$ appearing in this formula.

Let $R(x)$ be as in (\ref{SL2C matrix-function R}).
Formulae (\ref{GL2C transformation of the subprincipal symbol})
and
(\ref{definition of covariant subprincipal symbol})
imply that
the transformation (\ref{SL2C transformation of the operator})
of the differential operator induces the following transformation
of the matrix-function $L_\mathrm{csub}(x)$:
\begin{multline*}
L_\mathrm{csub}\mapsto
R^*(L_\mathrm{csub}+f(L_\mathrm{prin}))R
-f(R^*L_\mathrm{prin}R)
\\
+\frac i2
\left(
R^*_{x^\alpha}(L_\mathrm{prin})_{p_\alpha}R
-
R^*(L_\mathrm{prin})_{p_\alpha}R_{x^\alpha}
\right).
\end{multline*}
Comparing with
(\ref{SL2C transformation of the covariant subprincipal symbol})
we see that our function $f$ has to satisfy the condition
\begin{equation}
\label{condition on f}
f(R^*L_\mathrm{prin}R)
=
R^*f(L_\mathrm{prin})R
+\frac i2
\left(
R^*_{x^\alpha}(L_\mathrm{prin})_{p_\alpha}R
-
R^*(L_\mathrm{prin})_{p_\alpha}R_{x^\alpha}
\right)
\end{equation}
for any non-degenerate $2\times2$  Hermitian principal symbol
$L_\mathrm{prin}(x,p)$ and
any matrix-function (\ref{SL2C matrix-function R}).
Thus, we are looking for a function $f$ satisfying conditions
(\ref{homogeneity condition for f}) and (\ref{condition on f}).

Put
\begin{equation}
\label{formula for f}
f(L_\mathrm{prin}):=
-\frac i{16}\,
g_{\alpha\beta}
\{
L_\mathrm{prin}
,
\operatorname{adj}L_\mathrm{prin}
,
L_\mathrm{prin}
\}_{p_\alpha p_\beta},
\end{equation}
where subscripts $p_\alpha$, $p_\beta$ indicate partial derivatives and
\begin{equation}
\label{definition of generalised Poisson bracket}
\{F,G,H\}:=F_{x^\alpha}GH_{p_\alpha}-F_{p_\alpha}GH_{x^\alpha}
\end{equation}
is the generalised Poisson bracket on matrix-functions.
Note that the matrix-function in the RHS of formula (\ref{formula for f}) is Hermitian.

\begin{lemma}
\label{Lemma about f}
The function (\ref{formula for f})
satisfies conditions (\ref{homogeneity condition for f}) and (\ref{condition on f}).
\end{lemma}

\emph{Proof\ }
Substituting (\ref{psi transformation of the principal symbol})
into (\ref{formula for f})
we see that the terms with the gradient of the function $\psi(x)$
cancel out, which gives us (\ref{homogeneity condition for f}).
As to condition (\ref{condition on f}),
the appropriate calculations
are performed in
Appendix~\ref{Technical calculations I}.~$\square$

\

It is interesting that the generalised Poisson bracket on matrix-functions
(\ref{definition of generalised Poisson bracket})
was initially introduced for the purpose of abstract spectral analysis,
see formula (1.17) in~\cite{jst_part_a}.
It has now come handy in formula (\ref{formula for f}).

We will see later, in Section~\ref{Main result},
that the RHS of
(\ref{formula for f})
is just a way of writing the usual, Levi-Civita,
connection coefficients for spinor fields.
More precisely,
the RHS of
(\ref{formula for f})
does not give each spinor connection coefficient separately,
it rather gives their sum, the way they appear in the Dirac operator.

\begin{remark}
\label{uniqueness remark}
The function (\ref{formula for f}) is not a unique solution
of the system of equations
(\ref{homogeneity condition for f}) and (\ref{condition on f}):
one can always add $L_\mathrm{prin}(x,A(x))$, where $A(x)$ is
an arbitrary prescribed real-valued covector field.
We conjecture that our solution (\ref{formula for f})
of the system of equations
(\ref{homogeneity condition for f}) and (\ref{condition on f})
is unique up to the transformation
$f(L_\mathrm{prin})\mapsto
f(L_\mathrm{prin})+L_\mathrm{prin}(x,A(x))$.
Unfortunately, we are currently unable to provide a rigorous proof of this conjecture.
Moreover, even stating the uniqueness problem in a rigorous
and invariant fashion is a delicate issue.
Here the main difficulty is that our $f$ is not a function in the usual
sense, it is actually a nonlinear differential operator
mapping a $2\times2$ non-degenerate Hermitian principal symbol
$L_\mathrm{prin}(x,p)$
to a $2\times2$ Hermitian matrix-function
$(f(L_\mathrm{prin}))(x)$.
\end{remark}

\begin{remark}
If the conjecture stated in Remark~\ref{uniqueness remark} is true,
then the function (\ref{formula for f}) is singled out amongst all solutions
of the system of equations
(\ref{homogeneity condition for f}) and (\ref{condition on f})
by the property that it does not depend on any prescribed external fields.
\end{remark}

For the sake of clarity, we write down
the differential operator $L$ explicitly, in local coordinates,
it terms of its principal symbol $L_\mathrm{prin}$ and
covariant subprincipal symbol $L_\mathrm{csub}$.
Combining formulae
(\ref{operator in terms of principal and subprincipal symbols}),
(\ref{definition of covariant subprincipal symbol}) and
(\ref{formula for f}),
we get
\begin{multline}
\label{symbol Op explicit}
L=
-\frac i{2\sqrt{\rho(x)}}
\left(
[(L_\mathrm{prin})_{p_\alpha}(x)]\frac\partial{\partial x^\alpha}
+
\frac\partial{\partial x^\alpha}[(L_\mathrm{prin})_{p_\alpha}(x)]
\right)\sqrt{\rho(x)}
\\
-\frac i{16}
\left(
g_{\alpha\beta}
\{
L_\mathrm{prin}
,
\operatorname{adj}L_\mathrm{prin}
,
L_\mathrm{prin}
\}_{p_\alpha p_\beta}
\right)
\!(x)
+L_\mathrm{csub}(x).
\end{multline}
Here the covariant symmetric tensor $g_{\alpha\beta}(x)$ is the inverse
of the contravariant symmetric tensor $g^{\alpha\beta}(x)$ defined
by formula (\ref{definition of metric}),
$\{\,\cdot\,,\,\cdot\,,\,\cdot\,\}$ is the
generalised Poisson bracket on matrix-functions
defined by formula (\ref{definition of generalised Poisson bracket})
and $\,\operatorname{adj}\,$ is the operator of matrix adjugation
(\ref{definition of adjugation}).
See also Remark~\ref{correct reading} which explains how to
read formula (\ref{symbol Op explicit}) correctly.

\section{Electromagnetic covector potential}
\label{Electromagnetic covector potential}

The non-degeneracy condition (\ref{definition of non-degeneracy})
implies that for each $x\in M$ the matrices
$(L_\mathrm{prin})_{p_\alpha}(x)$, $\alpha=1,2,3,4$, form a basis
in the real vector space of $2\times2$ Hermitian matrices.
Here and throughout the paper the subscript $p_\alpha$ indicates partial differentiation.

Decomposing the covariant subprincipal symbol $L_\mathrm{csub}(x)$
with respect to this basis, we get
\begin{equation}
\label{decomposition of Z}
L_\mathrm{csub}(x)=(L_\mathrm{prin})_{p_\alpha}(x)\,A_\alpha(x)
\end{equation}
with some real coefficients $A_\alpha(x)$, $\alpha=1,2,3,4$.

Formula (\ref{decomposition of Z}) can be rewritten in more compact form as
\begin{equation}
\label{decomposition of Z compact}
L_\mathrm{csub}(x)=L_\mathrm{prin}(x,A(x)),
\end{equation}
where $A$ is a covector field with components $A_\alpha(x)$, $\alpha=1,2,3,4$.
Formula (\ref{decomposition of Z compact})
tells us that the covariant subprincipal symbol $L_\mathrm{csub}$
is equivalent to a real-valued covector field $A$,
the electromagnetic covector potential.

It is easy to see that our electromagnetic covector potential $A$ is invariant
under Lorentz transformations
(\ref{SL2C transformation of the operator})
and conformal scalings of the metric
(\ref{psi transformation of the operator}),
whereas formulae
(\ref{phi transformation of the covariant subprincipal symbol})
and
(\ref{decomposition of Z compact})
imply that
the transformation (\ref{phi transformation of the operator})
of the differential operator induces the transformation
\begin{equation}
\label{phi transformation of A}
A\mapsto A+\operatorname{grad}\phi.
\end{equation}

\section{Properties of the adjugate operator}
\label{Properties of the adjugate operator}

In this section we list gauge-theoretic properties of operator
adjugation~(\ref{definition of adjugate operator}).

Matrix adjugation
(\ref{definition of adjugation})
has the property
\begin{equation}
\label{group-theoretic property of matrix adjugation}
\operatorname{adj}(R^*PR)=R^{-1}(\operatorname{adj}P)(R^{-1})^*
\end{equation}
for any matrix $R\in\mathrm{SL}(2,\mathbb{C})$.
It is easy to see that operator adjugation
(\ref{definition of adjugate operator}) has a property
similar to (\ref{group-theoretic property of matrix adjugation}):
\begin{equation}
\label{group-theoretic property of operator adjugation}
\operatorname{Adj}(R^*LR)=R^{-1}(\operatorname{Adj}L)(R^{-1})^*
\end{equation}
for any matrix-function (\ref{SL2C matrix-function R}).

It is also easy to see that operator adjugation
(\ref{definition of adjugate operator}) commutes with the transformations
(\ref{psi transformation of the operator})
and
(\ref{phi transformation of the operator}):
\[
\operatorname{Adj}(e^\psi Le^\psi)=e^\psi(\operatorname{Adj}L)e^\psi,
\qquad
\operatorname{Adj}(e^{-i\phi}Le^{i\phi})=e^{-i\phi}(\operatorname{Adj}L)e^{i\phi}.
\]

Finally, let us observe that the map
(\ref{formula for f})
anticommutes with matrix adjugation
(\ref{definition of adjugation}),
\[
\operatorname{adj}f(L_\mathrm{prin})
=
-f(\operatorname{adj}L_\mathrm{prin}).
\]
This implies that the full symbol of the
operator $\operatorname{Adj}L$ is not necessarily the
matrix adjugate of the full symbol of the
operator $L$.

In the special case
when the principal symbol does not depend
on the position variable $x$ we get
$f(L_\mathrm{prin})=f(\operatorname{adj}L_\mathrm{prin})=0$,
so in this case
the full symbol of the
operator $\operatorname{Adj}L$ is the
matrix adjugate of the full symbol of the
operator $L$. The definition of the adjugate operator then simplifies
and becomes (\ref{definition of adjugate operator simplified}).

\section{Lorentz invariance of the operator (\ref{analytic definition of the Dirac operator})}
\label{Lorentz invariance}

In this section we show that our
Dirac operator (\ref{analytic definition of the Dirac operator})
is Lorentz invariant.
Recall that this operator acts on 4-columns of complex-valued scalar fields.

Let $R(x)$ be as in (\ref{SL2C matrix-function R}).
Define the $4\times4$ matrix-function
\begin{equation*}
\label{definition of matrix-function S}
S:=
\begin{pmatrix}
R&0\\
0&(R^{-1})^*
\end{pmatrix}.
\end{equation*}
Then
\begin{equation}
\label{analytic definition of the Dirac operator with tilde}
S^*DS=
\begin{pmatrix}
R^*LR&mI\\
mI&R^{-1}(\operatorname{Adj}L)(R^{-1})^*
\end{pmatrix}.
\end{equation}
The operator identity
(\ref{group-theoretic property of operator adjugation})
tells us that the diagonal terms in
(\ref{analytic definition of the Dirac operator with tilde})
are adjugates of each other, so formula
(\ref{analytic definition of the Dirac operator with tilde})
can be rewritten as
\begin{equation}
\label{analytic definition of the Dirac operator with tilde simplified}
S^*DS=
\begin{pmatrix}
R^*LR&mI\\
mI&\operatorname{Adj}(R^*LR)
\end{pmatrix}.
\end{equation}
We see that the operator
(\ref{analytic definition of the Dirac operator with tilde simplified})
has the same structure as
(\ref{analytic definition of the Dirac operator}),
which proves Lorentz invariance.

\section{Main result}
\label{Main result}

Formulae
(\ref{symbol Op explicit}),
(\ref{definition of Pauli matrices}),
(\ref{definition of metric}),
(\ref{definition of generalised Poisson bracket})
(\ref{decomposition of Z compact}),
(\ref{definition of adjugation})
and
(\ref{definition of adjugate operator})
allow us to rewrite our Dirac operator
(\ref{analytic definition of the Dirac operator})
in geometric notation --- in terms of Lorentzian metric,
Pauli matrices
and electromagnetic covector potential.
This raises the obvious question:
what is the relation between our Dirac operator
(\ref{analytic definition of the Dirac operator})
and the traditional Dirac operator
(\ref{traditional definition of the Dirac operator})?
The answer is given by the following theorem,
which is the main result of our paper.

\begin{theorem}
\label{main theorem}
Our Dirac operator
(\ref{analytic definition of the Dirac operator})
and the traditional Dirac operator
(\ref{traditional definition of the Dirac operator})
are related by the formula
\begin{equation}
\label{main theorem formula}
\rho^{1/2}
\,
D
\,
\rho^{-1/2}
=
|\det g_{\kappa\lambda}|^{1/4}
\,
D_\mathrm{trad}
\,
|\det g_{\mu\nu}|^{-1/4}\,,
\end{equation}
where $\rho$ is the density from our inner product
(\ref{definition of sesquilinear and bilinear forms}).
\end{theorem}

Here, of course, $\det g_{\kappa\lambda}=\det g_{\mu\nu}$.
We used different subscripts to avoid confusion because
tensor notation involves summation over repeated indices.

\

\emph{Proof of Theorem \ref{main theorem}\ }
Proving the $4\times4$ operator identity
(\ref{main theorem formula})
reduces to proving the following two separate $2\times2$
operator identities:
\begin{eqnarray}
\label{two by two operator identity 1}
\rho^{1/2}
\,
L
\,
\rho^{-1/2}
&\!\!\!=&\!\!\!
|\det g_{\kappa\lambda}|^{1/4}\,
\sigma^\alpha\,(-i\nabla+A)_\alpha\,|\det g_{\mu\nu}|^{-1/4}\,,
\\
\label{two by two operator identity 2}
\rho^{1/2}
\,
(\operatorname{Adj}L)
\,
\rho^{-1/2}
&\!\!\!=&\!\!\!
|\det g_{\kappa\lambda}|^{1/4}\,
\tilde\sigma^\alpha\,(-i\tilde\nabla+A)_\alpha\,|\det g_{\mu\nu}|^{-1/4}\,.
\end{eqnarray}
Here $\sigma^\alpha$ are Pauli matrices
(\ref{definition of Pauli matrices}),
$\tilde\sigma^\alpha$ are their adjugates,
and
$\nabla_\alpha$ and $\tilde\nabla_\alpha$
are covariant derivatives defined in accordance with formulae
(\ref{covariant derivative of undotted spinor field})
and
(\ref{covariant derivative of dotted spinor field}).

We shall prove the operator identity
(\ref{two by two operator identity 1}).
The operator identity
(\ref{two by two operator identity 2})
is proved in a similar fashion.

In the remainder of the proof we work in some local coordinate system.
The full symbols of the left-
and right-hand sides of
(\ref{two by two operator identity 1})
read
\[
(L_\mathrm{prin})_{p_\alpha}p_\alpha
-\frac i2(L_\mathrm{prin})_{x^\alpha p_\alpha}
-\frac i{16}\,
g_{\alpha\beta}
\{
L_\mathrm{prin}
,
\operatorname{adj}L_\mathrm{prin}
,
L_\mathrm{prin}
\}_{p_\alpha p_\beta}
+(L_\mathrm{prin})_{p_\alpha}A_\alpha
\]
and
\[
\sigma^\alpha p_\alpha
+\frac i4\sigma^\alpha(\ln|\det g_{\mu\nu}|)_{x^\alpha}
+\frac i4\sigma^\alpha\tilde\sigma_\beta
\left(
(\sigma^\beta)_{x^\alpha}
+\left\{{{\beta}\atop{\alpha\gamma}}\right\}\sigma^\gamma
\right)
+\sigma^\alpha A_\alpha
\]
respectively, where $\left\{{{\beta}\atop{\alpha\gamma}}\right\}$
denotes Christoffel symbols (\ref{definition of Christoffel symbols});
see also formulae
(\ref{operator in local coordinates})
and
(\ref{definition of the full symbol})
for the definition of the full symbol of a differential operator.
Comparing these with account of the fact that
$(L_\mathrm{prin})_{p_\alpha}=\sigma^\alpha$,
we see that the proof of the identity
(\ref{two by two operator identity 1})
reduces to the proof of the identity
\begin{multline}
\label{aux1}
-\frac i2(\sigma^\alpha)_{x^\alpha}
-\frac i{16}\,
g_{\alpha\beta}
\{
L_\mathrm{prin}
,
\operatorname{adj}L_\mathrm{prin}
,
L_\mathrm{prin}
\}_{p_\alpha p_\beta}
\\
=
\frac i4\sigma^\alpha(\ln|\det g_{\mu\nu}|)_{x^\alpha}
+\frac i4\sigma^\alpha\tilde\sigma_\beta
\left(
(\sigma^\beta)_{x^\alpha}
+\left\{{{\beta}\atop{\alpha\gamma}}\right\}\sigma^\gamma
\right).
\end{multline}
Using the standard formula
$(\ln|\det g_{\mu\nu}|)_{x^\alpha}=2\left\{{{\beta}\atop{\alpha\beta}}\right\}$
we rewrite (\ref{aux1}) as
\begin{multline}
\label{aux2}
\frac12\,
g_{\alpha\beta}
\{
L_\mathrm{prin}
,
\operatorname{adj}L_\mathrm{prin}
,
L_\mathrm{prin}
\}_{p_\alpha p_\beta}
\\
=
-2
\left(
2Ig^\alpha{}_\beta
+
\sigma^\alpha\tilde\sigma_\beta
\right)
\left(
(\sigma^\beta)_{x^\alpha}
+\left\{{{\beta}\atop{\alpha\gamma}}\right\}\sigma^\gamma
\right).
\end{multline}
Finally,
using formula (\ref{explicit formula for triple Poisson bracket})
we rewrite (\ref{aux2}) as
\begin{multline}
\label{agreement of covariant derivatives}
(\sigma^\alpha)_{x^\gamma}\tilde\sigma_\alpha\sigma^\gamma
-
\sigma^\gamma\tilde\sigma_\alpha(\sigma^\alpha)_{x^\gamma}
\\
=
-2
\left(
2Ig^\alpha{}_\beta
+
\sigma^\alpha\tilde\sigma_\beta
\right)
\left(
(\sigma^\beta)_{x^\alpha}
+\left\{{{\beta}\atop{\alpha\gamma}}\right\}\sigma^\gamma
\right).
\end{multline}

Thus, we have reduced the proof of the operator identity
(\ref{two by two operator identity 1})
to the proof of the identity
(\ref{agreement of covariant derivatives})
for Pauli matrices.
Calculations proving
(\ref{agreement of covariant derivatives})
are performed in
Appendix~\ref{Technical calculations II}.~$\square$

\

It remains only to note that
formula (\ref{main theorem formula}) implies
\begin{equation}
\label{main theorem formula rewritten}
D
=
\rho^{-1/2}
\,
|\det g_{\kappa\lambda}|^{1/4}
\,
D_\mathrm{trad}
\,
|\det g_{\mu\nu}|^{-1/4}\,
\rho^{1/2}\,.
\end{equation}
We identify a 4-column of complex-valued scalar fields $v$
with a bispinor field $\psi$
by means of the formula
\begin{equation}
\label{definition of bispinor}
v
=
|\det g_{\alpha\beta}|^{1/4}\,\rho^{-1/2}
\,\psi\,.
\end{equation}
Substituting
(\ref{main theorem formula rewritten})
and
(\ref{definition of bispinor})
into
(\ref{analytic statement of the Dirac equation})
we get
\begin{equation}
\label{analytic statement of the Dirac equation rewritten}
\rho^{-1/2}
\,
|\det g_{\kappa\lambda}|^{1/4}
\,
D_\mathrm{trad}
\,
\psi=0\,.
\end{equation}
Clearly, equation
(\ref{analytic statement of the Dirac equation rewritten})
is equivalent to equation
(\ref{Dirac equation traditional}).

\appendix

\section{Dirac equation in its traditional form}
\label{Dirac equation in its traditional form}

Before writing down the Dirac equation in its traditional form,
let us make several general  remarks on the notation that we will be using.
\begin{itemize}
\item
The notation in this appendix originates from \cite{LL4,buchbinder}.
Covariant derivatives of spinor fields are defined in accordance with
formulae (24) and (25) from \cite{JMP}.
The difference with \cite{LL4,buchbinder,JMP}
is that in the current paper
we enumerate local coordinates with indices $1,2,3,4$ rather than $0,1,2,3$.
Also, the difference with \cite{LL4,JMP}
is that in the current paper
we use opposite Lorentzian signature.
\item
The construction in this appendix is a generalisation of that
from Appendix A of \cite{jst_part_b}:
in \cite{jst_part_b} we dealt with the massless Dirac operator
in dimension three.
\item
We will write the Dirac equation in its \emph{spinor representation}
as opposed to its \emph{standard representation}, see Appendix B
in \cite{Mathematika} for details. The spinors $\xi^a$ and $\eta_{\dot b}$
that we will be using will be Weyl spinors,
i.e.~left-handed and right-handed spinors.
Let us note straight away that the $4\times4$ matrix differential
operator in the LHS of formula (B6) from \cite{Mathematika}
appears to have a structure different from
(\ref{analytic definition of the Dirac operator}). However,
it is easy to see that the representation
(B6) from \cite{Mathematika} reduces to
(\ref{analytic definition of the Dirac operator}) if one multiplies by
the constant $4\times4$ matrix $\begin{pmatrix}0&I\\ I&0\end{pmatrix}$
from the left.
\end{itemize}

The construction presented below is local, i.e.~we work
in a neighbourhood of a given point of a 4-manifold $M$ without boundary.
We have
a prescribed Lorentzian metric $g_{\alpha\beta}(x)$, $\alpha,\beta=1,2,3,4$,
and a prescribed electromagnetic covector potential $A_\alpha(x)$,
$\alpha=1,2,3,4$.  The metric tensor
is assumed to have three positive eigenvalues and one negative eigenvalue.

Consider a quartet of $2\times2$ Hermitian matrix-functions
$\sigma^\alpha{}_{\dot ab}(x)$.
Here the Greek index $\alpha=1,2,3,4$ enumerates the matrices, whereas
the Latin indices $\dot a=\dot 1,\dot 2$ and $b=1,2$ enumerate elements
of a matrix.
Here and throughout the appendix
the first spinor index always enumerates rows and the second columns.
We assume that under changes of local coordinates our
quartet of matrix-functions
transforms as the four components of a vector.
Throughout this appendix we use Greek letters for tensor indices
and we raise and lower tensor indices by means of the metric.

Define the ``metric spinor''
\begin{equation}
\label{metric spinor}
\epsilon_{ab}=\epsilon_{\dot a\dot b}=
\epsilon^{ab}=\epsilon^{\dot a\dot b}=
\begin{pmatrix}
0&-1\\
1&0
\end{pmatrix}.
\end{equation}
We will use the rank two spinor (\ref{metric spinor}) for raising and
lowering spinor indices. Namely, given
a quartet of $2\times2$ Hermitian matrix-functions
$\sigma^\alpha{}_{\dot ab}(x)$
we define the quartet of $2\times2$ Hermitian matrix-functions
$\tilde\sigma^{\alpha a\dot b}(x)$ as
\begin{equation}
\label{Hermitian matrix-functions with a tilde}
\tilde\sigma^{\alpha a\dot b}:=
-\epsilon^{ab}\,\epsilon^{\dot a\dot b}\,\sigma^\alpha{}_{\dot a b}\,.
\end{equation}
Note the order of spinor indices in the
matrix-functions $\tilde\sigma^{\alpha a\dot b}(x)$:
we choose it to be opposite to that in
\cite{JMP} but in agreement with that in
\cite{buchbinder}.

Examination of formulae
(\ref{metric spinor})
and
(\ref{Hermitian matrix-functions with a tilde})
shows that the  $2\times2$ matrices
$\sigma^\alpha{}_{\dot ab}$
and
$\tilde\sigma^{\alpha a\dot b}$
are adjugates of one another,
see formula (\ref{definition of adjugation}) for definition of matrix adjugation.
Hence, we could have avoided the use of the
``metric spinor'' in our construction of the Dirac equation,
using the mathematically more sensible concept of matrix adjugation instead.
The only reason we introduced the ``metric spinor''
is to relate the notation of the current paper to that of
\cite{LL4,buchbinder,JMP}.

Further on in this appendix we use matrix notation.
This means that we hide spinor indices and write the
matrix-functions
$\sigma^\alpha{}_{\dot ab}(x)$
and
$\tilde\sigma^{\alpha a\dot b}(x)$
as
$\sigma^\alpha(x)$
and
$\tilde\sigma^\alpha(x)$
respectively.

\begin{definition}
\label{abstract definition of Pauli matrices}
We say that the $2\times2$ Hermitian matrix-functions
$\sigma^\alpha(x)$ are \emph{Pauli matrices}
if these matrix-functions satisfy the identity
\begin{equation}
\label{defining identity for Pauli matrices}
\sigma^\alpha\tilde\sigma^\beta
+
\sigma^\beta\tilde\sigma^\alpha
=-2Ig^{\alpha\beta},
\end{equation}
where $I$ is the $2\times2$ identity matrix
and the tilde indicates matrix adjugation.
\end{definition}

\begin{remark}
The identity (\ref{defining identity for Pauli matrices}) is, of course, equivalent to
\begin{equation}
\label{defining identity for Pauli matrices 1}
\tilde\sigma^\alpha\sigma^\beta
+
\tilde\sigma^\beta\sigma^\alpha
=-2Ig^{\alpha\beta}.
\end{equation}
\end{remark}

Further on we assume that our $\sigma^\alpha(x)$ are Pauli matrices.

Consider a pair of spinor fields which we shall write as 2-columns,
\begin{equation}
\label{pair of spinor fields}
\xi=
\begin{pmatrix}
\xi^1\\
\xi^2
\end{pmatrix},
\qquad
\eta=
\begin{pmatrix}
\eta_{\dot1}\\
\eta_{\dot2}
\end{pmatrix}.
\end{equation}
Using matrix notation, we define the covariant derivatives
of these spinor fields as
\begin{equation}
\label{covariant derivative of undotted spinor field}
\nabla_\alpha\xi:=
\frac{\partial\xi}{\partial x^\alpha}
-\frac14\tilde\sigma_\beta
\left(
(\sigma^\beta)_{x^\alpha}
+\left\{{{\beta}\atop{\alpha\gamma}}\right\}\sigma^\gamma
\right)\xi\,,
\end{equation}
\begin{equation}
\label{covariant derivative of dotted spinor field}
\tilde\nabla_\alpha\eta:=
\frac{\partial\eta}{\partial x^\alpha}
-\frac14\sigma_\beta
\left(
(\tilde\sigma^\beta)_{x^\alpha}
+\left\{{{\beta}\atop{\alpha\gamma}}\right\}\tilde\sigma^\gamma
\right)\eta
\end{equation}
respectively, where
\begin{equation}
\label{definition of Christoffel symbols}
\left\{{{\beta}\atop{\alpha\gamma}}\right\}:=
\frac12g^{\beta\delta}
\left(
\frac{\partial g_{\gamma\delta}}{\partial x^\alpha}
+
\frac{\partial g_{\alpha\delta}}{\partial x^\gamma}
-
\frac{\partial g_{\alpha\gamma}}{\partial x^\delta}
\right)
\end{equation}
are the Christoffel symbols.

Formulae
(\ref{covariant derivative of undotted spinor field})
and
(\ref{covariant derivative of dotted spinor field})
warrant the following remarks.
\begin{itemize}
\item
The sign in front of the $\frac14$ in formula
(\ref{covariant derivative of undotted spinor field})
is the opposite of that in formula (24) of \cite{JMP}.
This is because in the current paper we use opposite
Lorentzian signature.
\item
The RHS of formula
(\ref{covariant derivative of undotted spinor field})
is a generalization of the expression appearing
in the RHS of formula (A.3) from \cite{jst_part_b}.
This follows from the observation that the adjugate
of a trace-free $2\times2$ matrix
$\sigma_\beta$ is $-\sigma_\beta$.
\item
If we multiply formula
(\ref{covariant derivative of undotted spinor field})
from the left by the ``metric spinor'' (\ref{metric spinor}),
apply complex conjugation and denote $\epsilon\bar\xi$ by $\eta$,
this gives us (\ref{covariant derivative of dotted spinor field}).
\end{itemize}

The massive Dirac equation reads
\begin{eqnarray}
\label{Dirac equation 1}
\sigma^\alpha\,(-i\nabla+A)_\alpha\,\xi+m\eta&\!\!\!=0\,,
\\
\label{Dirac equation 2}
\tilde\sigma^\alpha\,(-i\tilde\nabla+A)_\alpha\,\eta+m\xi&\!\!\!=0\,,
\end{eqnarray}
see formulae (B1) and (B2) from \cite{Mathematika}
or formulae (20.2) and (20.5) from \cite{LL4}.

We define the Dirac operator written in traditional geometric form as
\begin{equation}
\label{traditional definition of the Dirac operator}
D_\mathrm{trad}:=
\begin{pmatrix}
\sigma^\alpha\,(-i\nabla+A)_\alpha&mI\\
mI&\tilde\sigma^\alpha\,(-i\tilde\nabla+A)_\alpha
\end{pmatrix}
\end{equation}
and the bispinor field as the 4-column
\begin{equation}
\label{bispinor}
\psi:=
\begin{pmatrix}
\xi\\
\eta
\end{pmatrix}.
\end{equation}
Formulae
(\ref{Dirac equation 1})
and
(\ref{Dirac equation 2})
can then be rewritten as
\begin{equation}
\label{Dirac equation traditional}
D_\mathrm{trad}\,\psi=0\,.
\end{equation}

\section{Basic notions from the analysis of PDEs}
\label{Basic notions from the analysis of PDEs}

In this appendix we work with $m$-columns of complex-valued
scalar fields over an $n$-manifold $M$ without boundary. The main text of the
paper deals with the special case $n=4$, $m=2$, but in this
appendix $n$ and $m$ are arbitrary.

We assume that our manifold is equipped with a
prescribed positive density~$\rho$. This allows us to define
an inner product
on pairs $v$, $w$ of $m$-columns of
complex-valued scalar fields,
\begin{equation}
\label{definition of sesquilinear and bilinear forms}
\langle v,w\rangle:=\int_M w^*v\,\rho\,dx\,,
\end{equation}
where the star stands for Hermitian conjugation, $dx=dx^1\ldots dx^n$
and $x=(x^1,\ldots,x^n)$ are local coordinates.

Given a differential operator $L$, we define its formal adjoint $L^*$
by means of the formal identity
\begin{equation}
\label{definition of adjoint and transpose}
\langle Lv,w\rangle=\langle v,L^*w\rangle.
\end{equation}

Consider now a first order differential operator $L$.
In local coordinates it reads
\begin{equation}
\label{operator in local coordinates}
L=P^\alpha(x)\frac\partial{\partial x^\alpha}+Q(x),
\end{equation}
where $P^\alpha(x)$ and $Q(x)$ are some $m\times m$ matrix-functions
and summation is carried out over $\alpha=1,\ldots n$.
The full symbol of the operator $L$ is the matrix-function
\begin{equation}
\label{definition of the full symbol}
L(x,p):=iP^\alpha(x)\,p_\alpha+Q(x).
\end{equation}

Working with the full symbol is inconvenient
because the full symbol of a formally self-adjoint operator is not
necessarily Hermitian. The standard way of addressing this issue is
as follows.
We decompose
the full symbol into components homogeneous in $p$,
$\,L(x,p)=L_1(x,p)+L_0(x)\,$,
where
\begin{equation}
\label{definition of the homogeneous components of the full symbol}
L_1(x,p):=iP^\alpha(x)\,p_\alpha,
\qquad
L_0(x):=Q(x),
\end{equation}
and define the principal and subprincipal symbols as
\begin{equation}
\label{definition of the principal symbol}
L_\mathrm{prin}(x,p):=L_1(x,p),
\end{equation}
\begin{equation}
\label{definition of the subprincipal symbol}
L_\mathrm{sub}(x):=L_0(x)
+\frac i2(L_\mathrm{prin})_{x^\alpha p_\alpha}(x)
+\frac i2L_\mathrm{prin}(x,\operatorname{grad}(\ln\rho(x))),
\end{equation}
where $\rho$ is the density from
(\ref{definition of sesquilinear and bilinear forms}).
It is known that $L_\mathrm{prin}$ and $L_\mathrm{sub}$
are invariantly defined matrix-functions on $T^*M$ and $M$ respectively,
see subsection 2.1.3
in \cite{mybook} for details.

Let us explain why the formula for the subprincipal symbol
has the particular structure (\ref{definition of the subprincipal symbol}).
Firstly, using formulae (\ref{definition of the homogeneous components of the full symbol})
and (\ref{definition of the principal symbol})
we rewrite (\ref{definition of the subprincipal symbol}) as
\begin{equation}
\label{1}
L_\mathrm{sub}=Q
-\frac12(P^\alpha)_{x^\alpha}
-\frac12P^\alpha(\ln\rho)_{x^\alpha}\,.
\end{equation}
Here and further on in this paragraph
we drop, for the sake of brevity, the dependence on $x$.
The advantage of representing the subprincipal symbol in the form
(\ref{1}) is that the RHS is written explicitly in terms of the matrix-valued
coefficients $P^\alpha$ and $Q$ of the differential operator
(\ref{operator in local coordinates}).
Let us now substitute (\ref{operator in local coordinates}) into
the LHS of (\ref{definition of adjoint and transpose}),
use the formula for our inner product (\ref{definition of sesquilinear and bilinear forms})
and perform integration by parts.
We arrive at the expression for the adjoint operator in local coordinates
\begin{equation}
\label{2a}
L^*=\widehat P^\alpha\frac\partial{\partial x^\alpha}+\widehat Q,
\end{equation}
where
\begin{equation}
\label{2b}
\widehat P^\alpha
=-(P^\alpha)^*,
\qquad
\widehat Q=
Q^*
-[(P^\alpha)^*]_{x^\alpha}
-(P^\alpha)^*(\ln\rho)_{x^\alpha}\,.
\end{equation}
We then calculate the subprincipal symbol of $L^*$ using formula
(\ref{1}) and replacing matrix-valued
coefficients accordingly, compare formulae
(\ref{operator in local coordinates})
and
(\ref{2a}).
We get
\begin{equation}
\label{2c}
(L^*)_\mathrm{sub}=\widehat Q
-\frac12(\widehat P^\alpha)_{x^\alpha}
-\frac12\widehat P^\alpha(\ln\rho)_{x^\alpha}\,.
\end{equation}
Substitution of (\ref{2b}) into (\ref{2c}) gives us
\begin{equation}
\label{3}
(L^*)_\mathrm{sub}=
Q^*
-\frac12[(P^\alpha)^*]_{x^\alpha}
-\frac12(P^\alpha)^*(\ln\rho)_{x^\alpha}\,.
\end{equation}
Comparing formulae (\ref{1}) and (\ref{3}) we conclude that
\begin{equation}
\label{4}
(L^*)_\mathrm{sub}=(L_\mathrm{sub})^*.
\end{equation}
Thus, the whole point of introducing the two correction terms in
(\ref{definition of the subprincipal symbol}) (last two terms in the RHS)
is to ensure that we get the identity (\ref{4}).
Had we defined the subprincipal symbol as $L_\mathrm{sub}:=L_0$
we would not have the identity (\ref{4}).

The definition of the subprincipal symbol (\ref{definition of the subprincipal symbol})
originates from the classical paper \cite{DuiHor} of
J.J.~Duistermaat and L.~H\"ormander: see formula (5.2.8) in this paper.
Unlike \cite{DuiHor}, we work with matrix-valued symbols, but this
does not affect the formal definition of the subprincipal symbol.
What affects the definition of the subprincipal symbol is the fact that
we consider operators acting on columns of scalar fields rather
than operators acting on columns of half-densities and this leads 
to the appearance of the $\,\operatorname{grad}\ln\rho\,$ term in
(\ref{definition of the subprincipal symbol}).
Here we had to make a difficult decision: analysts prefer
to work with operators acting on half-densities because this
simplifies formulae, however the concept of a half-density is not
commonly used in the mathematical physics and theoretical physics communities.
We chose to avoid the use of the notion of a half-density
at the expense of having an extra correction term in 
(\ref{definition of the subprincipal symbol}).

For the principal symbol things are much easier and, obviously, we have
an analogue of formula  (\ref{4}):
\begin{equation}
\label{5}
(L^*)_\mathrm{prin}=(L_\mathrm{prin})^*.
\end{equation}

Examination of formulae
(\ref{operator in local coordinates})--(\ref{definition of the subprincipal symbol})
shows that
$L_\mathrm{prin}$, $L_\mathrm{sub}$ and $\rho$
uniquely determine the first order differential operator $L$.
Thus, the notions of principal symbol and subprincipal
symbol provide an invariant way of describing a first order
differential operator.

For the sake of clarity,
we write down the differential operator $L$ explicitly, in local coordinates,
in terms of its principal and subprincipal symbols:
\begin{multline}
\label{operator in terms of principal and subprincipal symbols}
L=
-\frac i{2\sqrt{\rho(x)}}
\left(
[(L_\mathrm{prin})_{p_\alpha}(x)]\frac\partial{\partial x^\alpha}
+
\frac\partial{\partial x^\alpha}[(L_\mathrm{prin})_{p_\alpha}(x)]
\right)\sqrt{\rho(x)}
\\
+L_\mathrm{sub}(x).
\end{multline}

\begin{remark}
\label{correct reading}
In writing formula
(\ref{operator in terms of principal and subprincipal symbols})
we used the convention that both operators of partial differentiation
$\frac\partial{\partial x^\alpha}$ act on all terms which come
(as a product) to the right, including the $m$-column of complex-valued
scalar fields $v$ which is present in
(\ref{operator in terms of principal and subprincipal symbols})
implicitly. Thus, a more explicit way of writing formula
(\ref{operator in terms of principal and subprincipal symbols})~is
\[
Lv=
-\frac{i(L_\mathrm{prin})_{p_\alpha}}{2\sqrt{\rho}}
\,
\frac{\partial(\sqrt{\rho}\,v)}{\partial x^\alpha}
-\frac{i}{2\sqrt{\rho}}
\,
\frac{\partial((L_\mathrm{prin})_{p_\alpha}\sqrt{\rho}\,v)}{\partial x^\alpha}
+L_\mathrm{sub}\,v\,.
\]
\end{remark}

\

Formulae (\ref{5}) and (\ref{4})
tell us that a first order differential operator is formally self-adjoint if and only
if its principal and subprincipal symbols are Hermitian matrix-functions.

We say that a formally self-adjoint first order differential operator $L$ is
\linebreak
\emph{elliptic} if
\begin{equation}
\label{definition of ellipticity}
\det L_\mathrm{prin}(x,p)\ne0,\qquad\forall(x,p)\in T^*M\setminus\{0\},
\end{equation}
and \emph{non-degenerate} if
\begin{equation}
\label{definition of non-degeneracy}
L_\mathrm{prin}(x,p)\ne0,\qquad\forall(x,p)\in T^*M\setminus\{0\}.
\end{equation}
The ellipticity condition (\ref{definition of ellipticity})
is a standard condition in the spectral theory of differential operators,
see, for example, \cite{jst_part_a}.
Our non-degeneracy condition
(\ref{definition of non-degeneracy})
is weaker and is designed to cover the case of hyperbolic operators.
In order to highlight the difference between
the ellipticity condition (\ref{definition of ellipticity})
and
the non-degeneracy condition (\ref{definition of non-degeneracy})
we consider two special cases.

\

\emph{Special case 1:}
$n=3$, $m=2$ and $\operatorname{tr}L_\mathrm{prin}(x,p)=0$.
In this case conditions
(\ref{definition of ellipticity})
and
(\ref{definition of non-degeneracy})
are equivalent.

\

\emph{Special case 2:}
 $n=4$ and $m=2$.
The proof of
Lemma~\ref{Lemma about Lorentzian metric}
shows that for each $x\in M$ there exists a
$p\in T_x^*M\setminus\{0\}$
such that $\det L_\mathrm{prin}(x,p)=0$, so it is impossible to satisfy
the ellipticity condition (\ref{definition of ellipticity}).
However, it is possible to satisfy
the non-degeneracy condition (\ref{definition of non-degeneracy}).
Indeed, consider the quantity
(density to the power $-1$)
$\,\det e_j{}^\alpha(x)$,
where $e_j$ is the frame from formula (\ref{principal symbol via frame 2}).
It is easy to see that
the non-degeneracy condition (\ref{definition of non-degeneracy})
is equivalent to the condition
$\det e_j{}^\alpha(x)\ne0$, $\forall x\in M$.
In other words,
the non-degeneracy condition (\ref{definition of non-degeneracy})
means that the vector fields $e_j$, $j=1,2,3,4$, encoded within
the principal symbol in accordance with formula
(\ref{principal symbol via frame 2})
are linearly independent at every point of our manifold $M$.

\section{Additional properties of Pauli matrices}
\label{Additional properties of Pauli matrices}

Throughout this appendix $\sigma^\alpha$, $\alpha=1,2,3,4$,
are Pauli matrices
and $\tilde\sigma^\alpha$ are their adjugates,
see Definition~\ref{abstract definition of Pauli matrices}.

\begin{lemma}
\label{Pauli lemma}
If $P$ is a $2\times2$ matrix then
\begin{equation}
\label{Pauli lemma formula 1}
\sigma_\alpha P\tilde\sigma^\alpha=
-2(\operatorname{tr}P)I,
\end{equation}
\begin{equation}
\label{Pauli lemma formula 2}
\sigma_\alpha P\sigma^\alpha=
2\operatorname{adj}P.
\end{equation}
\end{lemma}

\emph{Proof\ }
Formulae
(\ref{principal symbol via frame 2}),
(\ref{principal symbol via frame adjugate}),
(\ref{definition of Pauli matrices})
and
(\ref{definition of adjugate Pauli matrices})
imply
\begin{equation}
\label{Pauli lemma aux 1}
\sigma^\alpha=s^je_j{}^\alpha,
\qquad
\tilde\sigma^\alpha=s^j\tilde e_j{}^\alpha,
\end{equation}
where the matrices $s^j$ are defined in accordance with
(\ref{standard basis})
Substituting
(\ref{Pauli lemma aux 1})
into
(\ref{Pauli lemma formula 1})
and
(\ref{Pauli lemma formula 2})
and using the identities
(\ref{orthonormality of the frame})
and
(\ref{spatial inversion}),
we get
\[
\sigma_\alpha P\tilde\sigma^\alpha
=
-s^1Ps^1
-s^2Ps^2
-s^3Ps^3
-s^4Ps^4,
\]
\[
\sigma_\alpha P\sigma^\alpha
=
s^1Ps^1
+s^2Ps^2
+s^3Ps^3
-s^4Ps^4.
\]
The rest is a straightforward calculation.~$\square$

\

Note that an alternative way of proving formula
(\ref{Pauli lemma formula 1})
is by means of formula (1.2.27) from \cite{buchbinder}.

\section{Technical calculations I}
\label{Technical calculations I}

In this appendix we show that the function
(\ref{formula for f}) satisfies the condition (\ref{condition on f}).

Formulae
(\ref{definition of generalised Poisson bracket}),
(\ref{definition of Pauli matrices})
and
(\ref{definition of adjugate Pauli matrices})
give us
\begin{equation}
\label{explicit formula for triple Poisson bracket}
\frac12\,
g_{\alpha\beta}
\{
L_\mathrm{prin}
,
\operatorname{adj}L_\mathrm{prin}
,
L_\mathrm{prin}
\}_{p_\alpha p_\beta}
=
(\sigma^\alpha)_{x^\gamma}\tilde\sigma_\alpha\sigma^\gamma
-
\sigma^\gamma\tilde\sigma_\alpha(\sigma^\alpha)_{x^\gamma}.
\end{equation}
Note also that if we transform Pauli matrices $\sigma^\alpha$ as
\begin{equation}
\label{transformation of Pauli matrices}
\sigma^\alpha\mapsto R^*\sigma^\alpha R,
\end{equation}
where $R(x)$ is as in (\ref{SL2C matrix-function R}),
then the adjugate Pauli matrices $\tilde\sigma^\alpha$ transform~as
\begin{equation}
\label{transformation of adjugate Pauli matrices}
\tilde\sigma^\alpha\mapsto R^{-1}\tilde\sigma^\alpha(R^{-1})^*,
\end{equation}
see formula (\ref{group-theoretic property of matrix adjugation}).

Substituting formulae
(\ref{formula for f}),
(\ref{definition of Pauli matrices})
and
(\ref{explicit formula for triple Poisson bracket})--(\ref{transformation of adjugate Pauli matrices})
into
(\ref{condition on f})
we rewrite the latter as
$Q+Q^*=0$,
where
\begin{equation}
\label{Technical calculations I equation 1}
Q:=
-\frac i8
\left[
R^*\sigma^\alpha R_{x^\gamma}R^{-1}\tilde\sigma_\alpha\sigma^\gamma R
-
R^*\sigma^\gamma\tilde\sigma_\alpha\sigma^\alpha R_{x^\gamma}
\right]
+\frac i2R^*\sigma^\alpha R_{x^\alpha}.
\end{equation}
Hence, in order to prove (\ref{condition on f}) it is sufficient to prove
\begin{equation}
\label{Technical calculations I equation 1p1}
Q=0.
\end{equation}

Formula (\ref{defining identity for Pauli matrices 1}) implies that
$\tilde\sigma_\alpha\sigma^\alpha=-4I$,
so formula (\ref{Technical calculations I equation 1}) becomes
\begin{equation}
\label{Technical calculations I equation 2}
Q=
-\frac i8
R^*\sigma^\alpha R_{x^\gamma}R^{-1}\tilde\sigma_\alpha\sigma^\gamma R.
\end{equation}
The matrix-functions
$R_{x^\gamma}R^{-1}$ are trace-free,
so, by formula (\ref{Pauli lemma formula 1}),
\begin{equation}
\label{Technical calculations I equation 3}
\sigma^\alpha R_{x^\gamma}R^{-1}\tilde\sigma_\alpha=0.
\end{equation}
Formulae
(\ref{Technical calculations I equation 2})
and
(\ref{Technical calculations I equation 3})
imply
(\ref{Technical calculations I equation 1p1}).

\section{Technical calculations II}
\label{Technical calculations II}

In this appendix we prove the identity (\ref{agreement of covariant derivatives}).

Let us fix an arbitrary point $P\in M$ and prove the identity
(\ref{agreement of covariant derivatives}) at this point.
As the left- and right-hand sides
of (\ref{agreement of covariant derivatives}) are invariant under
changes of local coordinates~$x$,
it is sufficient to prove the identity
(\ref{agreement of covariant derivatives}) in
Riemann normal coordinates, i.e.~local coordinates
such that $x=0$ corresponds to the point $P$,
the metric at $x=0$ is Minkowski
and $\frac{\partial g_{\mu\nu}}{\partial x^\lambda}(0)=0$.
Moreover, as the identity we are proving involves only
first partial derivatives, we may assume, without loss of generality,
that the metric is Minkowski
for all $x$ in some neighbourhood of the origin.

Further on we assume that the metric is Minkowski.
We need to prove
\begin{equation}
\label{Technical calculations II eq 2}
Q=0\,,
\end{equation}
where
\begin{equation}
\label{Technical calculations II eq 3}
Q:=
(\sigma^\alpha)_{x^\gamma}\tilde\sigma_\alpha\sigma^\gamma
-
\sigma^\gamma\tilde\sigma_\alpha(\sigma^\alpha)_{x^\gamma}
+
2
\left(
2Ig^\alpha{}_\beta
+
\sigma^\alpha\tilde\sigma_\beta
\right)
(\sigma^\beta)_{x^\alpha}.
\end{equation}

Formula
(\ref{Technical calculations II eq 3})
can be rewritten in more compact symmetric form
\begin{equation}
\label{Technical calculations II eq 4}
Q=
(\sigma^\alpha)_{x^\gamma}\tilde\sigma_\alpha\sigma^\gamma
+
\sigma^\gamma\tilde\sigma_\alpha(\sigma^\alpha)_{x^\gamma}
+
4
(\sigma^\alpha)_{x^\alpha}.
\end{equation}
Using formulae
(\ref{defining identity for Pauli matrices}),
(\ref{defining identity for Pauli matrices 1})
and the fact that the metric is Minkowski
we can now rewrite
(\ref{Technical calculations II eq 4}) as
\begin{multline}
\label{Technical calculations II eq 5}
Q=
(\sigma^\alpha)_{x^\gamma}
(
-2g_\alpha{}^\gamma
-\tilde\sigma^\gamma\sigma_\alpha
)
+
(
-2g^\gamma{}_\alpha
-\sigma_\alpha\tilde\sigma^\gamma
)
(\sigma^\alpha)_{x^\gamma}
+
4
(\sigma^\alpha)_{x^\alpha}
\\
=
-(\sigma^\alpha)_{x^\gamma}\tilde\sigma^\gamma\sigma_\alpha
-\sigma_\alpha\tilde\sigma^\gamma(\sigma^\alpha)_{x^\gamma}
=
-(\sigma^\alpha)_{x^\gamma}\tilde\sigma^\gamma\sigma_\alpha
-\sigma^\alpha\tilde\sigma^\gamma(\sigma_\alpha)_{x^\gamma}
\\
=
\sigma^\alpha
(\tilde\sigma^\gamma)_{x^\gamma}
\sigma_\alpha
-
(
\sigma^\alpha
\tilde\sigma^\gamma
\sigma_\alpha
)_{x^\gamma}.
\end{multline}
Formula (\ref{Pauli lemma formula 2}) allows us to rewrite formula
(\ref{Technical calculations II eq 5}) in the form
\[
Q=
2
\left[
\operatorname{adj}
\left(
(\tilde\sigma^\gamma)_{x^\gamma}
\right)
-
(\operatorname{adj}\tilde\sigma^\gamma)_{x^\gamma}
\right].
\]
As the operations of matrix adjugation
(\ref{definition of adjugation})
and partial differentiation commute, we arrive at
(\ref{Technical calculations II eq 2}).

\end{document}